\documentclass[12pt]{amsart}
\usepackage{amsmath}
\usepackage{amsthm}
\usepackage{amssymb}
\usepackage{comment}
\newcommand{\n}{\par\noindent}

\newcommand{\sn}{\par\smallskip\noindent}
\newcommand{\mn}{\par\medskip\noindent}
\newcommand{\bn}{\par\bigskip\noindent}
\newcommand{\pars}{\par\smallskip}
\newcommand{\parm}{\par\medskip}

\newtheorem{theorem}{Theorem}

\newtheorem{proposition}[theorem]{Proposition}
\newtheorem{corollary}[theorem]{Corollary}

\newtheorem{remar}[theorem]{Remark}

\newcommand{\gzaun}{\unskip\nobreak\hfil\penalty50%
\hskip1em\hbox{}\nobreak\hfil%
$\#$\parfillskip=0pt\finalhyphendemerits=0}

\newcommand{\ac}{^{\rm ac}}
\newcommand{\chara}{\mbox{\rm char}\,}
\newcommand{\trdeg}{\mbox{\rm trdeg}\,}
\newcommand{\Gal}{\mbox{\rm Gal}\,}

\newcommand{\isom}{\simeq}
\newcommand{\bfind}[1]{\index{#1}{\bf #1}}

\newcommand{\cO}{\mathcal{O}}
\newcommand{\cM}{\mathcal{M}}

\newcommand{\cE}{\mathcal{E}}

\newcommand{\cC}{\mathcal{C}}

\newcommand{\cL}{\mathcal{L}}

\newcommand{\Z}{\mathbb Z}
\newcommand{\N}{\mathbb N}
\newcommand{\Q}{\mathbb Q}
\newcommand{\R}{\mathbb R}
\newcommand{\F}{\mathbb F}
\begin{document}            

\title[Model theory of tame valued fields and beyond]{Model theory of tame valued fields
and beyond: recent developments and open questions}
\author{Franz-Viktor Kuhlmann}
\address{Institute of Mathematics, University of Szczecin,
ul. Wielkopolska 15,
70-451 Szczecin, Poland}
\email{fvk@usz.edu.pl}

\date{28.\ 11.\ 2025}

\begin{abstract}
We give a survey on recent developments in the model theory of valued fields since
the introduction of the notion of ``tame valued field'', and of the modifications and
generalizations of this notion.
\end{abstract}

\thanks{I thank Konstantinos Kartas, Arno Fehm, Franziska Jahnke and Blaise Boissonneau
for helpful
discussions and input to this survey. I thank the Institute for Mathematical Science
of the National University of Singapore, where I gave the talk on which this survey
paper is based, for support and hospitality. Further, I thank the Hausdorff Research
Institute for Mathematics in Bonn for support and hospitality, as my stay at the
Institute during the program ``Definability, decidability, and computability'' gave
me the opportunity to discuss the contents of this paper with several colleagues.}


\maketitle

%
%
\section{Introduction}
The notion of ``tame valued field'' was first introduced in the author's PhD thesis;
for an extended version of it, see \cite{Kuth2}. The name was coined in collaboration
with the author's very supportive supervisor Peter Roquette. Since then, tame valued
fields and their model theory have been generalized, modified or applied in many
research articles. One of the applications was to prove local uniformization of
Abhyankar places in positive characteristic (\cite{Ku23}), and of local uniformization
by alteration of arbitrary places in positive characteristic (\cite{Ku26}). However,
the developments in the model theory of valued fields have been much more complex, and
they are what this survey is devoted to. It is based on a survey talk given at the
conference {\it Recent Applications of Model Theory} held in June of 2025 at the
Institute for Mathematical Science of the National University of Singapore.
For the slides of the talk, see
{\tt https://www.fvkuhlmann.de/TalkSingapore2025slides.pdf}, 
and
{\tt https://www.fvkuhlmann.de/AbstractSingapore2025.pdf}
for the  extended abstract.

\parm
Through the course of this survey, the reader will meet tame, separably tame and roughly
tame fields, extremal fields, perfectoid fields, deeply ramified and roughly deeply
ramified fields. Without any claim for completeness, several main theorems and
contributors are presented and 16 open problems are listed. I hope that this will
provide a useful basis for young as well as experienced mathematicians eager to attack
these problems.

%
%
\section{Some preliminaries}
For a valued field $(K,v)$, we denote its value group by $vK$,
its residue field by $Kv$, and its valuation ring by $\cO_K\,$.
By $(L|K,v)$ we denote an extension $L|K$ with a valuation $v$ on $L$, where $K$
is endowed with the restriction of $v$. In this case, there are induced embeddings of
$vK$ in $vL$ and of $Kv$ in $Lv$. The extension $(L|K,v)$ is called \bfind{immediate}
if these embeddings are onto. A valued field $(K,v)$ is called \bfind{algebraically
maximal} if it does not admit nontrivial immediate algebraic extensions, and it is
called \bfind{maximal} if it does not admit any nontrivial immediate extensions.

A valued field $(K,v)$ is called \bfind{henselian} if for each algebraic extension
$L|K$ the extension of $v$ to $L$ is unique.
A finite extension $(L|K,v)$ of a henselian valued field $(K,v)$ is called
\bfind{defectless} if
\[
[L:K]\>=\>(vL:vK)[Lv:Kv]\>.
\]
This definition can be extended to valued fields that are not henselian by using
the \bfind{fundamental inequality} (cf.\ (17.5) of \cite{End} or Theorem 19 on p.~55
of \cite{ZS}): given a valued field $(K,v)$ and a finite field extension $L|K$,
there are finitely many extensions of $v$ to $L$ and we have
\begin{equation}                             \label{fundineq}
[L:K]\>\geq\>\sum_{i=1}^{\rm g} (v_i L:vK)[Lv_i:Kv]\>,
\end{equation}
where $v_1,\ldots,v_{\rm g}$ are the distinct extensions of $v$ from $K$ to $L$.
We say that $(K,v)$ is \bfind{defectless in $L$} if equality holds in (\ref{fundineq}).
A valued field is called a \bfind{defectless field} if it is defectless in each finite extension. A valued field is a defectless field if and only if any (and then every)
of its henselizations is (\cite[Theorem 8.9]{KuHab}).
\pars
The following are examples for defectless fields:
\sn
1) \ every trivially valued and every algebraically closed valued field,
\n
2) \ every valued field with residue characteristic 0,
\n
3) \ maximal fields (see the discussion at the beginning of Section 4 in \cite{Ku58});
in particular:
\n
-- with its $p$-adic valuation $v_p\,$, the field $\Q_p$ of $p$-adic numbers,
\n
-- with their $t$-adic valuations $v_t\,$, the Laurent series field $k((t))$ over a
field $k$ and every power series field $k((t^\Gamma))$ with coefficients in $k$ and
exponents in $\Gamma$ (see \cite[Theorem 8.26 and Corollaries 8.27 and 8.28]{KuHab}),
\n
4) \ $(\Q,v_p)$ and its henselization (by \cite[Theorem 8.32]{KuHab}),
\n
5) \ $(k(t),v_t)$ and its henselization (this is a special case of \cite[Theorem 1.1]{Ku29}).

%
%
\section{Tame fields}
An algebraic extension $(L|K,v)$ of a henselian valued field $(K,v)$ is called
\bfind{tame} if every finite subextension $K'|K$ satisfies the following conditions:
\sn
(T1) \ the ramification index $(vK':vK)$ is not divisible by $\chara Kv$,
\sn
(T2) \ the residue field extension $K'v|Kv$ is separable,
\sn
(T3) \ the extension $(K'|K,v)$ is defectless.

\sn
{\bf Remark.}
This notion of ``tame extension'' does not coincide with the notion of
``tamely ramified extension'' as defined on page 180 of O.~Endler's book
\cite{End}. The latter definition requires (T1) and
(T2), but not (T3). Our tame extensions are the defectless tamely
ramified extensions in the sense of Endler's book. In particular, in our
terminology, proper immediate algebraic extensions of henselian fields
are not tame (in fact, they cause a lot of problems in the model
theory of valued fields).

\pars
A henselian valued field $(K,v)$ is called a \bfind{tame field} if the algebraic
closure $K\ac$ of $K$ with the unique extension of $v$ is a tame extension of $(K,v)$.
It follows from conditions (T1)--(T3) that all tame fields are
perfect. For a characterization of tame fields, see Theorem~\ref{tf} below.

\pars
The perfect hull $\F_p((t))^{1/p^\infty}$ of the field $\F_p((t))$ of formal Laurent
series over the field $\F_p$ with $p$ elements is perfect but not tame, as the
extension generated by a root of the Artin-Schreier polynomial $X^p-X-1/t$ is of
degree $p$ and immediate (see \cite[Example 3.12]{Ku31}), and therefore does not
satisfy (T3).

For details on tame fields, see \cite{Ku39}. The results of that paper are now
frequently applied in the model theory of valued fields. In particular:
\begin{theorem}[Kuhlmann (2016)]          \label{MTtf}
Tame fields $(K,v)$ satisfy model completeness  in the language $\cL_{\rm val}$ of
valued rings  relative to the elementary theories of their value groups $vK$ in the
language of ordered groups  and their residue fields $Kv$ in the language of rings.
If $\chara K=\chara Kv$, then also relative completeness  and relative decidability
hold.
\end{theorem}

However, there are still daunting questions about tame fields  that have remained
unanswered.

\sn
{\bf Open problem 1:} Do tame fields admit quantifier elimination in a suitable
language?
\sn
The problem is that we do not know enough about  \bfind{purely wild extensions}, i.e.,
algebraic extensions of a henselian valued field that are linearly disjoint
from tame extensions. For background, see \cite{Ku2}.

\parm
There are also open problems about the model theory of tame fields $(K,v)$  that have
\bfind{mixed characteristic}, i.e., $\chara K=0$ while $\chara Kv=p>0$. One question is
whether (or under which additional conditions) they satisfy relative completeness and
decidability. In the article \cite{Ku46}, examples are given of two tame fields
$(K_1,v_1)$ and $(K_2,v_2)$  with $v_1K_1\equiv v_1K_2$ and $K_1v_1\equiv K_2v_2$
such that $(K_1,v_1)\not\equiv (K_2,v_2)$. The difference to the case of tame fields
of positive characteristic is that in this case the restriction of the
valuations to the prime fields are trivial  while in the mixed characteristic case they
are $p$-adic and hence nontrivial.

Progress on this problem has been made in \cite[Theorem 1.2]{L}. Lisinski proves:
\begin{theorem}[Lisinski (2021)]             \label{ThmLi1}
Take tame fields $(L,v)$ and $(F,w)$ of mixed characteristic with residue characteristic
$p>0$ such that $vL\equiv wF$ in the language of ordered groups  with a constant symbol
$\pi$ interpreted as $v(p)$ and $w(p)$, respectively,  and that $Lv\equiv Fw$ in the
language of rings.  Assume that the relative algebraic closure $(K,v)$ of $\Q$ in $(L,v)$
is algebraically maximal  and that $vL/vK$ is torsion free.  Assume further that every
monic polynomial $f\in \Z[X]$  has a root in $\cO_F$ if it has a root in $\cO_L\,$.
Then $(L,v)\equiv (F,w)$ in $\cL_{\rm val}\,$.
\end{theorem}

As mentioned already in Theorem~\ref{MTtf}, relative decidability is
proven for tame fields of positive characteristic. As an application, we obtain:
\begin{theorem}[Kuhlmann (2016)]
Take $q=p^n$ for some prime $p$ and some $n\in\N$, and an ordered abelian group
$\Gamma$.  Assume that $\Gamma$ is divisible  or elementarily equivalent to the
$p$-divisible hull of $\Z$.  Then the $\cL_{\rm val}$-elementary theory  of the
power series field $\F_q((t^\Gamma))$ with coefficients in the field $\F_q$ with $q$
elements and exponents in $\Gamma$, endowed with its canonical valuation $v_t\,$, is
decidable.
\end{theorem}

Lisinski improves this result as follows (\cite[Theorem 1]{L}):
\begin{theorem}[Lisinski (2021)]             \label{ThmLi2}
Take a perfect field $\F$ of characteristic $p>0$  whose elementary theory in the
language of rings is decidable, and a $p$-divisible ordered abelian group $\Gamma$
whose elementary theory in the language of ordered groups  with a constant symbol $1$
is decidable. Then the $\cL_{\rm val}(t)$-elementary theory of $\F((t^\Gamma))$ is
decidable.
\end{theorem}
\n
Here, $\cL_{\rm val}(t)$ denotes the language $\cL_{\rm val}$ with a constant
symbol~$t$.
\pars
Lisinski also proves a theorem giving a criterion for two tame fields containing
$\F_p(t)$  to be equivalent in $\cL_{\rm val}(t)$ that is analogous to
Theorem~\ref{ThmLi1} (see \cite[Theorem 1.1]{L}).

%
%
\section{The fields $\Q_p$, $\F_p((t))$ and their algebraic extensions}
Since in 1965 Ax and Kochen in \cite{AK}, and independently Ershov in \cite{E},
established the decidability of the elementary theory of $\Q_p$, several questions
about the decidability of the elementary or the existential theory of local fields
and their extensions have been answered,  and several others have remained
open. We have already seen some results in equal positive charateristic. In
contrast, less is known in mixed characteristic, for instance about
\n
$\bullet$ \ the totally ramified extension $\Q_p(\zeta_{p^\infty})$ obtained from
$\Q_p$   by adjoining all $p^n$-th roots of unity, $n\in\N$,
\n
$\bullet$ \ the totally ramified extension $\Q_p(p^{1/p^\infty})$ obtained from
$\Q_p$  by adjoining a compatible system of $p^n$-th roots of $p$, $n\in\N$,
\n
$\bullet$ \ the maximal abelian extension $\Q_p^{ab}$ of $\Q_p\,$.

\sn
These fields with their canonical valuations are studied in \cite{Ka2}.
\begin{theorem}[Kartas (2024)]    \label{Kart1}
The fields $\Q_p(\zeta_{p^\infty})$ and $\Q_p(p^{1/p^\infty})$ equipped with their unique
extensions $v_p$ of the $p$-adic valuation  admit maximal immediate extensions which
have decidable elementary $\cL_{\rm val}$-theories.
\end{theorem}
\n
It was shown by W.~Krull in \cite{Kr} that every valued field $(K,v)$ admits a
maximal immediate extension $(M,v)$ (the proof was later simplified by
K.~A.~H.~Gravett in \cite{Gra}).
All of these maximal immediate extensions are tame fields.  But the fields themselves
are not \bfind{Kaplansky fields}, i.e., fields satisfying hypothesis (A) in \cite{Ka},
and Kartas shows that there are uncountably many maximal
immediate  extensions with distinct elementary $\cL_{\rm val}$-theories.  This
implies that uncountably many of them are not decidable.

Kartas proves a ``perfectoid transfer theorem''  (\cite[Theorem A]{Ka2}) which transfers
the decidability in certain expansions of $\cL_{\rm val}$ of fields in equal positive
characteristic to the decidability in $\cL_{\rm val}$ of suitable untilts.

By Theorem~\ref{ThmLi2},  the perfectoid field $\F_p((t^\Gamma))$,
where $\Gamma$ is the $p$-divisible hull of $\Z$,  is decidable in the language
$\cL_{\rm val}(t)$. Kartas constructs a suitable untilt $K$ of $\F_p((t^\Gamma))$
which by the perfectoid transfer theorem is decidable in the language
$\cL_{\rm val}\,$. As $\F_p((t^\Gamma))$ is a maximal immediate extension  of the
completion of the perfect hull $\F_p(t^{1/p^\infty})$, which is the tilt of the
completion of $\Q_p(p^{1/p^\infty})$, a theorem of Fargues and Fontaine can be used to
show that $K$  is a maximal immediate extension of the latter  and hence also of
$\Q_p(p^{1/p^\infty})$ itself. The case of $\Q_p(\zeta_{p^\infty})$ is
similar.  This completes the proof of Theorem~\ref{Kart1}.

\pars
Kartas notes that all tilts of the undecidable maximal immediate extensions of $\Q_p
(\zeta_{p^\infty})$ and $\Q_p(p^{1/p^\infty})$  are maximal immediate extensions of
$\F_p((t))^{1/p^\infty}$. Being tame fields, they are decidable in the language
$\cL_{\rm val}\,$.  But they are not decidable in the language $\cL_{\rm val}(t)$.
\sn
{\bf Open problem 2:} What is the structure of these extensions (apart from the fact that
they are infinite)?  What are the indications in their structure  that distinguish
the decidable from the undecidable extensions?
\pars
Kartas also notes that $\Q_p^{ab}$, being a Kaplansky field, admits a unique maximal
immediate extension,
and that it follows from the model theory of algebraically maximal Kaplansky fields
that this extension is decidable in $\cL_{\rm val}\,$.

\parm
{\bf Open problem 3:} Are $\Q_p(\zeta_{p^\infty})$, $\Q_p(p^{1/p^\infty})$ and
$\Q_p^{ab}$  decidable in $\cL_{\rm val}\,$?  Are $\F_p((t))^{1/p^\infty}$ and
$\F_p^{\rm ac}((t))^{1/p^\infty}$  decidable in $\cL_{\rm val}$ or even
$\cL_{\rm val}(t)$?
\sn
Here $\F_p^{\rm ac}$ denotes the algebraic closure of $\F_p\,$.

\pars
We do not know the answers,  but there are some conditional
results connecting decidability of the mixed characteristic fields  with those of the
positive characteristic fields.  Although these fields are not perfectoid,  Kartas
succeeds to deduce the following  from the perfectoid transfer theorem.
\begin{theorem}[Kartas (2024)]
\mbox{ }\n
(a) If $\,\F_p((t))^{1/p^\infty}$ has a decidable elementary or existential
$\cL_{\rm val}(t)$-theory, then $\Q_p(\zeta_{p^\infty})$ and $\Q_p(p^{1/p^\infty})$
have decidable elementary or existential $\cL_{\rm val}$-theories, respectively.
\sn
(b) If $\,\F_p^{\rm ac}((t))^{1/p^\infty}$ has a decidable elementary or existential
$\cL_{\rm val}(t)$-theory, then $\Q_p^{ab}$ has a decidable elementary or existential
$\cL_{\rm val}$-theory, respectively.
\end{theorem}

\sn
{\bf Open problem 4:} What about the reverse direction?
\mn
In fact, if $\Q_p(\zeta_{p^\infty})$ or $\Q_p(p^{1/p^\infty})$ has a decidable
$\cL_{\rm val}$-theory, then $\F_p((t))^{1/p^\infty}$ has a decidable
$\cL_{\rm val}$-theory, and if $\Q_p^{ab}$ has a decidable $\cL_{\rm val}$-theory,
then so does $\F_p^{\rm ac}((t))^{1/p^\infty}$. This essentially follows from
\cite[Corollary 1.7.6]{JK}, which says that if a perfectoid field $K$ has a decidable
$\cL_{\rm val}$-theory, then so does its tilt. Since problem 4 is about non-complete
valued fields, one also needs to use that the involved fields have the same elementary
theory as their respective completions. So a more precise version of problem 4 is whether
this argument can be improved to get decidability in $\cL_{\rm val}(t)$ on the positive
characteristic side. (For a related result, compare \cite[Proposition 7.2.3]{JK}.)

\parm
Let us point out that $\Q_p(\zeta_{p^\infty})$, $\Q_p(p^{1/p^\infty})$ and $\Q_p^{ab}$
with their unique extensions of the $p$-adic valuation are not tame extensions of
$\Q_p\,$. This brings us to the question: what can be said about infinite tame extensions
of $\Q_p\,$? In \cite{Ka2}, Kartas works with a precise formulation of resolution of
singularities, which he calls ``Log-Resolution''. By the work of Hironaka, it is known
that Log-Resolution holds in characteristic 0. Under the assumption that it also holds
in positive characteristic, Kartas proves an existential Ax-Kochen/Ershov principle
(\cite[Theorem A]{Ka1}), from which the following decidability results for infinite tame extensions of $\Q_p$ and $\F_p((t))$ can be deduced:
\sn
$\bullet$ \ for primes $\ell\ne p$, \ $\Q_p(p^{1/\ell^\infty})$ has a decidable
existential theory in the language of rings,
\sn
$\bullet$ \ for primes $\ell\ne p$, \ $\F_p((t))(t^{1/\ell^\infty})$ has a decidable
existential theory in the language of rings enriched by a constant symbol~$t$,
\sn
$\bullet$ \ the maximal tame extensions (also known as absolute ramification fields,
(see \cite[Proposition 4.1]{Ku2}) of $\Q_p$ and $\F_p((t))$ have decidable existential
theories in the respective languages.

%
%
\section{Separably tame fields}
A henselian field is called a \bfind{separably tame field}
if every separable-algebraic extension is a tame extension.  We let  $\cL_{\rm val, Q}$
denote the language $\cL_{\rm val}$  enriched by $m$-ary predicates $Q_m\,$, $m\in\N$,
for $p$-independence.  That is, in a field $K$ of characteristic $p>0$,  $Q_m$ is
interpreted by
\[
Q_m(x_1, \ldots, x_m)\>\Leftrightarrow\>\left\{
\begin{array}{l}
\mbox{the monomials of exponents $< p$ in the $x_i$'s}\\
\mbox{are linearly independent over the subfield $K^p$}\\
\mbox{of $p$-th powers.}
\end{array}
\right.
\]

\sn
Field extensions $L|K$ as $\cL_{\rm val, Q}$-structures are separable,  i.e., linearly
disjoint from the perfect hull of $K$.

The model theory of separably tame fields is studied in the article \cite{Ku45}.
\begin{theorem}[Kuhlmann -- Pal (2016)]                  \label{ThmKP}
Separably tame fields $(K,v)$ of positive characteristic and finite degree of
inseparability  satisfy completeness and decidability in $\cL_{\rm val}$
relative to the elementary theories of their value groups $vK$ in the language of ordered
groups  and of their residue fields $Kv$ in the language of rings.

In the language $\cL_{\rm val, Q}\,$, they also satisfy relative model completeness.
\end{theorem}

In \cite[Corollary 1.6]{A}, Anscombe removes the condition of finite degree of
inseparability from the relative decidability result in $\cL_{\rm val}$, and in
\cite[Theorem 1.5]{A} from the other assertions of Theorem~\ref{ThmKP}
in a language $\cL_{{\rm val,}\lambda}$
which is $\cL_{\rm val,Q}$  with the predicates $Q_m$
replaced by function symbols for Lambda functions (see \cite[Definition 2.5]{A}).
To this end, Anscombe proves that the $\cL_{{\rm val,}\lambda}$-theory of equal
characteristic separably tame valued fields has the \bfind{Lambda Relative Embedding
Property} (see \cite[Definition 4.10]{A}). This is done by adapting the
proofs of \cite[Theorem 7.1]{Ku39}, which shows that the elementary class
of tame fields has the \bfind{Relative Embedding Property} (see \cite[Section 6]{Ku39}),
and of \cite[Teorem 5.1]{Ku45}, which shows that the elementary class of separably
tame fields of finite degree of inseparability has the \bfind{Separable Relative
Embedding Property} (see \cite[Section 4]{Ku45}).

%
%
\section{Perfectoid and deeply ramified fields}
In \cite[Definition 1.2]{Sch} Peter Scholze defines a \bfind{perfectoid field} to be a
complete nondiscrete rank 1 valued field of residue characteristic $p>0$ such that
$\cO_K/p\cO_K$ is \bfind{semiperfect}, that is, the Frobenius is surjective on
$\cO_K/p\cO_K\,$. This implies that the value group is $p$-divisible. A valued field
has \bfind{rank 1} if its value group is embeddable in the ordered additive group $\R$.

Neither ``complete'' nor ``rank 1'' are elementary properties.
A suitable elementary class of valued fields containing the perfectoid fields
is that of deeply ramified fields, studied in the article \cite{Ku68}.
A nontrivially valued field $(K,v)$ is a \bfind{deeply ramified field} if and only
if the following conditions hold:
\sn
{\bf (DRvg)} whenever $\Gamma_1\subsetneq\Gamma_2$ are convex subgroups of the value
group $vK$,  then $\Gamma_2/\Gamma_1$ is not isomorphic to $\Z$  (that is, no
archimedean component of $vK$ is discrete),
\sn
{\bf (DRvr)} if $\chara Kv=p>0$, then $\cO_K/p\cO_K$ is semiperfect if $\chara K=0$, and
$\cO_{\widehat K}/p\cO_{\widehat K}$ is semiperfect if $\chara K=p$, where
$\cO_{\widehat K}$ is the valuation ring of the completion $\widehat K$ of $(K,v)$.
\sn
If $(K,v)$ has rank 1, then (DRvg) just means that $(K,v)$ is not discrete. If $(K,v)$
is complete,  then (DRvr) means that $\cO_K/p\cO_K$ is semiperfect. Hence every
perfectoid field is deeply ramified. Every perfect valued field of positive
characteristic $p$ (in particular, $\F_p((t))^{1/p^\infty}$) and every tame field is
deeply ramified. The former as well as all tame fields of residue characteristic $p>0$
have $p$-divisible value group, but this does not
necessarily hold for deeply ramified fields of characteristic $0$ with residue
characteristic $p>0$. In \cite{Ku68}, we define $(K,v)$ to be a \bfind{semitame
field} if it is a deeply ramified field whose value group is $p$-divisible if $\chara
Kv=p>0$. Semitame fields form a smaller elementary class which still contains all
perfectoid fields.

%
%
\section{Roughly deeply ramified and roughly tame fields}
Inspired by the notion ``roughly $p$-divisible value group'' introduced by Will Johnson,
we call $(K,v)$ a \bfind{roughly deeply ramified field},  or in short
an \bfind{rdr field},  if it satisfies axiom (DRvr) together with:
\sn
{\bf (DRvp)} if $\chara Kv=p>0$,  then $v(p)$ is not the smallest positive element
in the value group $vK$.
\sn
The two axioms (DRvp) and (DRvr) together  imply that the smallest convex subgroup of
$vK$  containing $v(p)$  (or equivalently, the interval $[-v(p),v(p)]\,$)  is
$p$-divisible.

\parm
For the definition of roughly tame fields, we need the following characterization of
tame fields as preparation. The following is part of \cite[Theorem 3.2]{Ku39}:
\begin{theorem}[Kuhlmann (2016)]    \label{tf}
A henselian field $(K,v)$ is a tame field if and only if the following conditions
hold:
\sn
(TF1) \ if $\chara Kv=p>0$, then $vK$ is $p$-divisible,
\sn
(TF2) \ $Kv$ is perfect,
\sn
(TF3) \ $(K,v)$ is algebraically maximal.
\end{theorem}
\sn
Replacing (TF1) by
\sn
{\it (TF1r) \ if $\chara Kv=p>0$, then $[-v(p),v(p)]$ is $p$-divisible,}
\sn
we obtain the definition of a \bfind{roughly tame field}.

\pars
In the article \cite{RS}, the following is proven:
\begin{theorem}[Rzepka -- Szewczyk (2023)]         \label{thmRS}
A henselian field is roughly tame  if and only if all of its algebraic extension
fields  are defectless fields.
\end{theorem}

We know that being defectless is an important property  in the model theory of valued
fields.  For instance, if a valued field is existentially closed in its maximal
immediate extensions,  then it is henselian and defectless (see \cite[Lemma 5.5]{Ku39}).
The converse is not true, see Proposition~\ref{hdfneciie}.

The property of being a defectless field is preserved under finite algebraic extensions,
but in general not under infinite algebraic extensions. For instance, Example 3.12 of
\cite{Ku31} constructs infinite algebraic extensions of $(\F_p(t),v_t)$ and
$(\F_p((t)),v_t)$ which are not defectless
fields, and Example 3.20 of \cite{Ku31} constructs infinite algebraic extensions of
$(\Q,v_p)$ and $(\Q_p,v_p)$ which are not defectless fields.

%
%
\section{Taming perfectoid fields}
In the article \cite{JK}, Jahnke and Kartas generalize the model theoretic results
about tame fields to the elementary class of roughly tame fields. To this end, they
prove the Relative Embedding Property for the elementary class of roughly tame fields
(\cite[Fact 3.3.12]{JK}). From this it follows that the assertions of Theorem~\ref{MTtf}
also hold for the elementary class of roughly tame fields. They put this generalization
to work in their approach  of ``taming perfectoid fields''.  They work with an
elementary class $\cC$ of henselian fields $(K,v)$  of residue characteristic $p>0$
with distinguished element  $\pi\in K\setminus\{0\}$, $v\pi>0$, such that:
\sn
$(\cC 1)$ \ the valuation ring $\cO_K$ is semitame, and
\n
$(\cC 2)$ \ with the coarsening $w$ of $v$ associated with the valuation ring $\cO_v
[\pi^{-1}]$,  $(K,w)$ is algebraically maximal  (which implies that it is roughly
tame).
\sn
(See \cite[Definition 4.2.1 and Proposition 4.2.2]{JK}.)
The class $\cC$ contains all henselian roughly deeply ramified fields of mixed
characteristic (see \cite[Remark 4.2.4]{JK}).

\pars
If $(K',v')$ is a suitable ultrapower of a perfectoid field $(K,v)$ with distinguished
element  $\varpi\in K\setminus\{0\}$, $v\varpi>0$, and $w'$ is the coarsest coarsening
of $v'$ on $K'$ such that $w'\varpi>0$, then $(K',w')\in\cC$ for any $\pi\in K
\setminus\{0\}$, $v\pi>0$, and the residue field $K'w'$ with its valuation induced by
the one of the ultrapower is an elementary extension of the tilt of $K$
(see \cite[Corollary 4.2.6 and Theorem 6.2.3]{JK}). This is used to show that certain
model theoretic properties hold for perfectoid fields if and only if they hold for
their tilts (see \cite[Corollary 5.3.1]{JK}).
\pars
For the class $\cC$, Jahnke and Kartas prove analogues  of the model theoretic
results for (roughly) tame fields, but with the residue fields $Kv$ replaced  by
the residue rings $\cO_K/\pi\cO_K$ (see \cite[Theorems 5.1.2 and 5.1.4]{JK}). This
``mods out the non-tame part'' of the valued fields in $\cC$. So we are still left
with the
\sn
{\bf Open problem 5:} What can we say about the model theory of (roughly) deeply
ramified fields (relative to their value groups and residue fields), and in
particular of $\F_p((t))^{1/p^\infty}$?

\pars
It is well known that the henselization $(\F_p(t)^h,v_t)$ of $(\F_p(t),v_t)$ is
existentially closed in $(\F_p((t)),v_t)$. This can be deduced from the following more
general result (see \cite[Theorem 5.12]{Ku39}):
\begin{theorem}[Kuhlmann 2016]
Let $k$ be an arbitrary field. Then $(k(t),v_t)^h$ is existentially closed in
$(k((t)),v_t)$.
\end{theorem}

However, the following has remained a daunting
\sn
{\bf Open problem 6:} Is $\F_p(t)^h$ an elementary substructure of $\F_p((t))$?
\sn
In contrast, Jahnke and Kartas prove that $\F_p(t^{1/p^\infty})^h$  is an elementary
substructure of $\F_p((t))^{1/p^\infty}$ (\cite[Corollary 1.7.5]{JK}). This positive
result encourages us to ask:
\sn
{\bf Open problem 7:} Is it possible to prove model theoretic results  for henselian
perfect valued fields of positive characteristic,  analoguous to those for tame
fields  (but under mild additional conditions)?

%
%
\section{On the model theory of $\F_p((t))$}
While model theoretic results about $\Q_p$  and in particular the decidability of
$\Q_p$ are known  since the work of Ax--Kochen and Ershov,  we are still facing the
\sn
{\bf Open problem 8:} What can we say about the model theory  and in particular
a complete axiomatization and the decidability  of $\F_p((t))$?
\pars
In the article \cite{Ku13} the following negative result is proven:
\begin{theorem}[Kuhlmann (2001)]
The $\cL_{\rm val}(t)$-elementary axiom system
\sn
(A$_t$) ``henselian defectless valued field of positive characteristic  with value
group a $\Z$-group with smallest element $v(t)$  and residue field~$\F_p$''
\sn
is not complete.
\end{theorem}
\n
This theorem is proven by constructing an extension $(L,v)$ of $(\F_p((t)),v_t)$
with the following properties:
\sn
$\bullet$ \ $(L,v)$ satisfies axiom system (A$_t$),
\sn
$\bullet$ \ $L|K$ is of trancendence degree 1 and regular (i.e., $L|K$ is separable
and $K$ is relatively algebraically closed in $L$),
\sn
$\bullet$ \ there is an $\forall\exists$-elementary $\cL_{\rm val}(t)$-sentence
expressing a property of additive polynomials which holds in $(K,v)$ but not in $(L,v)$,
\sn
see \cite[Theorem 1.3]{Ku13}.
A polynomial $f(X)\in K[X]$ is called \bfind{additive} if $f(a+b)=f(a)+f(b)$ for all
$a,b$ in any extension field of $K$. If $\chara K=p>0$, then the additive polynomials
in $K[X]$ are precisely the polynomials of the form
\[
\sum_{i=0}^{m} c_i X^{p^i}\;\;\;\mbox{ with } c_i\in K\,,\, m\in\N
\]
(see \cite[VIII, \S 11]{L}). If $K$ is infinite, then $f(X)\in K[X]$ is additive if and
only if $f(a+b)=f(a)+f(b)$ for all $a,b\in K$. Additive polynomials in several variables
are defined in a similar way; but note that they are just sums of additive polynomials
in one variable. The special role of
additive polynomials for valued fields of characteristic $p>0$ had been long known; for
example, the \bfind{Artin-Schreier polynomial} $X^p-X$ is additive.

\pars
Moreover, it is shown that $(L,v)$ is not $\cL_{\rm val}$-existentially closed in its
maximal immediate extensions (cf.\ \cite[Theorem 1.3]{Ku13}). This proves:
\begin{proposition}[Kuhlmann (2001)]              \label{hdfneciie}
There are henselian defectless fields that are not $\cL_{\rm val}$-existentially
closed in their maximal immediate extensions.
\end{proposition}
\sn
{\bf Open problem 9:} Is there a handy additional condition on the immediate
extensions that remedies this situation?

\pars
In the article \cite{Ku13}, also the following is shown:
\begin{theorem}[Kuhlmann (2001)]
The $\cL_{\rm val}$-elementary axiom system
\n
(A) ``henselian defectless valued field of positive characteristic  with value
group a $\Z$-group  and residue field~$\F_p$''
\n
is not complete.
\end{theorem}
\n
Further, a (not really handy) axiom scheme, called (PDOA) and stating properties of
additive polynomials, is suggested to be added to axiom systems (A) or (A$_t$).
A much more elegant axiom scheme was found
after Yuri Ershov introduced the notion of ``extremal field''  and claimed that
$\F_p((t))$ is extremal.  However, his definition and proof were faulty.  In the
article \cite{Ku35} it is shown that $\F_p((t))$ does not satisfy Ershov's definition;
in fact, every valued field satisfying this definition must be algebraically closed.
A corrected definition is given,  and it is shown that $\F_p((t))$ satisfies this
corrected definition,  which we present now.

%
%
\section{Extremal fields}
A valued field $(K,v)$ is called \bfind{extremal} if for every
multi-variable polynomial $f(X_1, \dots , X_n)$ over $K$,
the set
\[
\{v(f(a_1,\dots,a_n))\mid a_1,\dots, a_n \in \cO_K\}
\subseteq vK\cup\{\infty\}
\]
has a maximal element.  This is an $\cL_{\rm val}$-elementary axiom scheme. Ershov's error was to put ``$K$'' in place of ``$\cO_K$''.
\begin{theorem}[Azgin -- Kuhlmann -- Pop (2012)]
$\F_p((t))$ is an extremal field.
\end{theorem}
\sn
{\bf Open problem 10:} Is (A) + ``$(K,v)$ is extremal'' a complete axiom system?

\parm
In the article \cite{Ku46} an almost complete characterization of extremal valued field
is given:
\begin{theorem}[Anscombe -- Kuhlmann (2016)]
Let $(K,v)$ be a nontrivially valued field. If $(K,v)$ is extremal, then
it is henselian and defectless, and
\begin{itemize}
  \item[(i)] $vK$ is a $\mathbb{Z}$-group, or
  \item[(ii)] $vK$ is divisible and $Kv$ is large.
 \end{itemize}
Conversely, if $(K,v)$ is henselian and defectless, and
\begin{itemize}
  \item[(i)] $vK\isom\Z$, or $vK$ is a $\mathbb{Z}$-group and $\chara Kv=0$, or
  \item[(ii)] $vK$ is divisible and $Kv$ is large and perfect,
 \end{itemize}
then $(K,v)$ is extremal.
\end{theorem}

\sn
A complete characterization is not known, and there are many more open problems about
extremal fields listed in \cite{Ku46}. For instance, in contrast to properties such as
``henselian'' and ``defectless'', we do not entirely know how extremality behaves under
composition of valuations:
\sn
{\bf Open problem 11:} If $v=w\circ\bar{w}$ with $w$ and $\bar{w}$ extremal and $w$
has divisible value group, does it follow that $v$ is extremal?

\sn
{\bf Open problem 12:} We know that if $v=w\circ\bar{w}$ is extremal, then so is
$\bar{w}$ (see \cite[Lemma 4.1]{Ku46}). But does it also follow that $w$ is extremal?

\pars
On the other hand, it is shown in \cite[Theorem 1.13]{Ku46} that in a certain sense,
extremal fields are abundant in valuation theory:
\begin{theorem}[Anscombe -- Kuhlmann (2016)]
Let $(K,v)$ be any $\aleph_{1}$-saturated valued field. Assume that $\Gamma$ and $\Delta$
are convex subgroups of $vK$ such that $\Delta\subsetneq\Gamma$ and $\Gamma/\Delta$ is
archimedean. Let $u$ (respectively $w$) be the coarsening of $v$ corresponding to
$\Delta$ (resp. $\Gamma$). Denote by $\bar{u}$ the valuation induced on $Kw$ by
$u$. Then $(Kw,\bar{u})$ is maximal, extremal and large, and its value group is
isomorphic either to $\Z$ or to $\R$. In the latter case, also  $Ku=(Kw)\bar{u}$ is
large.
\end{theorem}
This shows that extremal fields can be seen as the rank 1 building blocks of valuations,
at least of those which are $\aleph_{1}$-saturated. The properties of the valuations
thus built (which can vary wildly) apparently depend on the way the building blocks are
``glued together'', a process that remains to be investigated further.

\pars
Finally, let us mention that in the extension $(L,v)$ of $\F_p((t))$ constructed in
\cite{Ku13} not all images of additive polynomials have the optimal approximation
property (for its definition, see \cite[\S3]{Ku46}). It thus follows from
\cite[Theorem 3.4]{Ku46} that $(L,v)$ is not extremal. This gives rise to the following
\sn
{\bf Open problem 13:} Is every extremal valued field existentially closed in its maximal immediate extensions?

%
%
\section{The existential theory of $\F_p((t))$ and criteria for large fields to be
existentially closed in extensions}
Let us return to the model theory of $\F_p((t))$. The following is shown in
\cite{AF}:
\begin{theorem}[Anscombe -- Fehm (2016)]       \label{ThmAF}
The existential $\cL_{\rm val}$-theory of $\F_p((t))$ is decidable.
\end{theorem}
\n
However, as can be seen from our discussion of Kartas' work, we would like to have
more. In the article \cite{DS}, Jan Denef and Hans Schoutens proved in 2003 that the
existential $\cL_{\rm val}(t)$-theory of $\F_p((t))$ is decidable, provided that
resolution of
singularities holds in positive characteristic.  In order to discuss more recent
improvements of this result,  we need some preparations.

\parm
In the article \cite{Ku20}
the following question is studied:  Take a field extension $F|K$  such that $F$
admits a \bfind{$K$-rational place},  or in other words, a valuation with residue field
$K$.  Under which additional conditions does it follow  that $K$ is existentially
closed in $F$? Here a key role is played by large fields.  While they are usually
defined in a different way (see \cite[Section 1.3]{Ku20}), one can also use the
model theoretic approach: A field $K$ is \bfind{large} if it is existentially closed
in $K((t))$.
\begin{theorem}[Kuhlmann (2004)]         \label{R4perf}
Let $K$ be a perfect field.  Then the following conditions are equivalent:
\n
1) \ $K$ is a large field,
\n
2) \ $K$ is existentially closed in every power series field $K((t^\Gamma))$,
\n
3) \ $K$ is existentially closed in every extension field $L$ which
admits a $K$-rational place.
\end{theorem}

\pars
What can we say about fields that are not perfect? For a first answer, we introduce a
hypothesis that enables us to generalize the above theorem.
\bfind{Local uniformization} is a local form, and a consequence, of resolution of
singularities. For background, see \cite{Ku23, Ku26, T}. In recent decades, doubts
have spread in the community of algebraic geometers working on resolution of
singularities that it can be proven for all dimensions in positive characteristic.
However, there is much more hope for a corresponding general version of local
uniformization.
\begin{theorem}[Kuhlmann (2004)]
If all rational places of arbitrary function fields admit local uniformization,
then the three conditions of Theorem~\ref{R4perf} are equivalent, for arbitrary
fields~$K$.
\end{theorem}

\parm
In the paper \cite{ADF}
the assumption that implication 1)$\Rightarrow$3) holds for {\it arbitrary} fields
$K$  is called hypothesis \bfind{(R4)}. Hence local uniformization implies (R4). The
authors prove the following strengthening of Theorem~\ref{ThmAF}:
\n
\begin{theorem}[Anscombe -- Dittmann -- Fehm (2023)]
If (R4) holds, then the existential $\cL_{\rm val}(t)$-theory of $\F_p((t))$ is
decidable.
\end{theorem}
\sn
{\bf Open problem 14:} Does (R4) hold?

\parm
In the paper \cite{dJ} de Jong proved resolution by \bfind{alteration}, which means that
a finite extension of the function field  of the algebraic variety under consideration
is taken into the bargain.  By valuation theoretical tools, Knaf and Kuhlmann proved
local uniformization by alteration in \cite{Ku26}.
\sn
{\bf Open problem 15:} Does local uniformization by alteration imply  a reasonable
(and useful) hypothesis ``(R4) by alteration''?  Is there a ``model theory by
alteration''?
\mn
You cannot always get what you want -- but
perhaps after a finite extension?
\sn
Indeed, the following is shown in \cite{KK3}:
\begin{theorem}[Kuhlmann -- Knaf (?)]
Take a large field $K$ and a function field $F|K$ which admits a rational place. Then
there is a finite purely inseparable extension $K'|K$ such that $K'$ is existentially
closed in the field compositum $F.K'$.
\end{theorem}
\n
One possible proof uses Temkin's ``inseparable local uniformization'' by alteration
(\cite[Theorem 1.3.2]{T}. However, Arno Fehm has pointed out that the theorem can also
more directly be deduced from results in \cite{Ku20}.

\parm
In \cite{CP}, Cossart and Piltant prove resolution of singularities for algebraic
varieties of dimension no larger than $3$. This implies that locl uniformization holds for all places on function fields of transcendence degree no larger than $3$. In this situation, we can employ \cite[Theorem 19]{Ku20}:
\begin{theorem}                             \label{kecF1}
Let $K$ be a large field and $F|K$ an algebraic function field. If there
is a rational place of $F|K$ which admits local uniformization, then $K$
is existentially closed in~$F$.
\end{theorem}

From this, we obtain:
\begin{corollary}
Let $K$ be a large field and $F|K$ an algebraic function field with $\trdeg F|K\leq
3$. If there is a rational place of $F|K$ which admits local uniformization, then $K$
is existentially closed in~$F$.
\end{corollary}

Likewise, we obtain that if $F|K$ is an algebraic function field over a large field
$K$ admitting a rational place, then $K$ is existentially closed in~$F$ with respect
to all existential sentences with at most three existential quantifiers.

%
%
\section{Classification of defects}
Let me give some details on the classification of
defects which has been introduced in \cite{Ku68}. If $(L|K,v)$ is a finite extension
for which the extension of $v$ from $K$ to $L$ is unique, then by the Lemma of Ostrowski
(\cite[Corollary to Theorem 25, Section G, p.\ 78]{ZS}),
\begin{equation}                    \label{feuniq}
[L:K]\>=\> \tilde{p}^{\nu }\cdot(vL:vK)[Lv:Kv]\>,
\end{equation}
where $\nu$ is a non-negative integer and $\tilde{p}$ the
\bfind{characteristic exponent} of $Kv$, that is, $\tilde{p}=\chara Kv$ if it is
positive and $\tilde{p}=1$ otherwise. The factor $d(L|K,v):=\tilde{p}^{\nu }$ is
the \bfind{defect} of the extension $(L|K,v)$. By our previous definition, if
$(K,v)$ is henselian, then $(L|K,v)$ is
a defectless extension if $d(L|K,v)=1$. This always holds if $\chara Kv=0$.

\pars
Take a Galois extension $\cE=(L|K,v)$ of prime degree $p$ with nontrivial defect; then
$p=\chara Kv$. For every $\sigma$ in its Galois group $\Gal (L|K)$, with
$\sigma\ne\,$id, we set
\begin{equation}                        \label{Sigsig}
\Sigma_\sigma\>:=\> \left\{ v\left( \left.\frac{\sigma f-f}{f}\right) \right| \,
f\in L^{\times} \right\} \>.
\end{equation}
This set is a final segment of $vK$ and independent of the choice of $\sigma$
(see \cite[Theorems 3.4 and 3.5]{Ku68}); we denote it by $\Sigma_\cE\,$.
We say that $\cE$ has \bfind{independent defect} if
\begin{equation}                                     \label{indepdef}
\left\{\begin{array}{lcr}
\Sigma_{\cE}\!\!&=&\!\!\! \{\alpha\in vK\mid \alpha >H_\cE\}\>\mbox{ for some proper
convex subgroup $H_\cE$}\\
&&\!\!\! \mbox{ of $vK$ such that $vK/H_\cE$ has no smallest positive element;}
\end{array}\right.
\end{equation}
otherwise we say that $\cE$ has \bfind{dependent defect}. If $(K,v)$ has rank 1, then
condition (\ref{indepdef}) just means that $\Sigma_{\cE}$ consists of all positive
elements in $vK$.

If $(K,v)$ is of mixed characteristic, then we set $K':=K(\zeta_p)$, where $\zeta_p$ is
a primitive $p$-th root of unity; otherwise, we set $K':=K$. Note that every Galois
extension $L$ of prime degree $p$ of a field $K$ of characteristic 0 containing a
primitive $p$-th root of unity is a \bfind{Kummer extension}, i.e., it is generated
by an element $\eta$ with $\eta^p\in K$. Now we call $(K,v)$ an \bfind{independent
defect field} if for some extension of $v$ to $K'$, all Galois extensions of $(K',v)$
of degree $p$  with nontrivial defect have independent defect. This definition does
not depend on the chosen extension of $v$ as all extensions are conjugate.
\pars
The following is Theorem 1.5 of \cite{Ku68}:
\begin{theorem}[Kuhlmann -- Rzepka 2023]                \label{algext}
Every algebraic extension of a deeply ramified field is again a deeply ramified field.
The same holds for semitame fields and for roughly deeply ramified fields.
\end{theorem}
\n
Further, every roughly deeply ramified field is an independent defect field (see
\cite[Theorem 1.10 (1)]{Ku68}). This proves:
\begin{theorem}[Kuhlmann -- Rzepka 2023]
Every algebraic extension of a roughly deeply ramified field is an independent defect
field.
\end{theorem}
\n
In view of Theorem~\ref{thmRS}, we conjecture that a henselian field is a roughly deeply
ramified field if and only if all of its algebraic extensions are independent defect
fields. Thus we ask:
\sn
{\bf Open problem 16:} If all algebraic extensions of a henselian field $(K,v)$
are independent defect fields, does it follow that $(K,v)$ is a roughly deeply
ramified field?

%
%
\section{Some definable valuation rings}
It is obvious from its definition that a Galois extension $(L|K,v)$ of prime degree
with independent defect
gives rise to a coarsening $\cO_{\cE,K}$ of the valuation ring $\cO_K\,$, namely
\[
\cO_{\cE,K}\>:=\>\{b\in K\mid \exists\alpha\in H_\cE:\>\alpha\leq vb\}
\]
whose value group is $vK/H_\cE\,$. It is shown in \cite{Kudef} that $\cO_{\cE,K}$
is $\cL_{\rm val}$-definable when $L|K$ is an Artin-Schreier or a Kummer extension.
For this it is not needed that $(K,v)$ be henselian, and in fact, it is applied
to deeply ramified fields, which are not required to be henselian. For the case of
henselian $(K,v)$, the above is used in \cite[Theorem 4.11]{KRS} to define
corresponding henselian valuations on $K$ that are already definable in the language
of rings. It follows that a henselian roughly deeply ramified field which is not
defectless always has a henselian valuation definable in the language
of rings (cf.\ \cite[Corollary 4.14]{KRS}).

\pars
For ur own applications, however, we are more interested in definable coarsening
$\cO_{\cE}$ of the valuation ring $\cO_L\,$. In the case of Galois extensions
$(L|K,v)$ of prime degree that either have nontrivial defect or satisfy $[L:K]=
(vL:vK)$, corresponding coarsenings $\cO_{\cE}$ can be defined in the language
$\cL_{{\rm val},K}$ of valued fields with a predicate for membership in $K$. These
coarsenings, together with their maximal ideals $\cM_{\cE}\,$, play an important role
in our study of Galois defect extensions of prime degree (see \cite{Kudef}). For
instance, in \cite[Theorem 1.4]{pr1}, $\cM_\cE$ is used, under the notation
$\cM_{v_\cE}\,$, in the characterization of extensions with independent defect. In
\cite{pr2} it is used to present the K\"ahler differential $\Omega_{\cO_L|\cO_K}$ of
the extension.

\bn\bn

\end{document}